\documentclass[10pt]{article}
\usepackage{amssymb}
\usepackage{enumerate}
\usepackage{graphicx}

\def\R{{\rm I\! R}}

\def\psx{{\partial p \over \partial x}}
\def\qsx{{\partial q \over \partial x}}
\def\rsx{{\partial r \over \partial x}}
\def\psy{{\partial p \over \partial y}}
\def\qsy{{\partial q \over \partial y}}
\def\rsy{{\partial r \over \partial y}}

\newtheorem{theorem}{Theorem}
\newtheorem{corollary}{Corollary}

\newtheorem{lemma}{Lemma}

\newtheorem{remark}{Remark}

\title{Real divergence-free Jacobian maps are shear maps  
\thanks{
2000 AMS {\it Mathematics Subject Classification}: 14R15.
{\it Key words and phrases:}  Jacobian Conjecture, divergence-free, shear map.}
}

\author{M. Sabatini}
\begin{document}
\maketitle
\begin{abstract} We show that divergence-free jacobian maps, recently considered in \cite{N}, \cite{S},  are actually shear maps.\end{abstract}

\section{Introduction}

Let $\Phi:\R^n \rightarrow \R^n$ be a continuously differentiable map. By the inverse function theorem, if its Jacobian determinant $\det J_\Phi$ does not vanish at a point, then $\Phi$ is locally invertible at such a point. It is well-known that the local invertibility of $\Phi$, even if it holds on all of $\R^n$, does not imply $\Phi$'s global invertibility. Additional conditions implying $\Phi$'s global invertibility, or just injectivity, have been studied in several fields (see \cite{AP}, \cite{P}). Some conditions have just been conjectured to be sufficient to guarantee the map's invertibility, as in the case of the celebrated Jacobian Conjecture \cite{K}. Originally formulated for complex maps, and successively studied for maps defined on arbitrary  fields, such a conjecture asks to prove that a polynomial map with constant non-zero Jacobian determinant is globally invertible, with polynomial inverse. 
The problem has been attacked in several ways. 
One can find in \cite{BCW} an overview of results concerning the Jacobian Conjecture  up to 1982. A more recent account is contained in \cite{dBvdE}.  Among general results concerning such a problem, it is known that it is equivalent to prove or disprove the statement in any field of zero characteristic, that it is sufficient to prove $\Phi$'s injectivity in order to get its surjectivity \cite{BCW}, and that $\Phi$'s global invertibility implies that $\Phi^{-1}$ is a polynomial map. 
The most studied special case is the bidimensional one, $\Phi (x,y) = (P(x,y) ,Q(x,y) )$, where the statement was proved under the hypothesis that either P's or Q's degree is 4, or prime, or both degrees are $\leq$ 100 (see \cite{BCW} for a more comprehensive list of results). 

Recently, a particular class of maps has been taken into account, namely the so-called divergence-free maps \cite{N}. 
Such maps have the form $\Phi (x,y) = (x+p(x,y),y+q(x,y))$, where $p(x,y)$ and $q(x,y)$ are polynomials of order $> 1$, with $\psx + \qsy =0$. For such maps one has $\det J_\Phi = 1 + \left(\psx + \qsy \right) + \left(\psx \qsy - \psy \qsx \right) $, hence under the condition $\psx + \qsy =0$ also the determinant-like  term  vanishes, $\psx \qsy - \psy \qsx =0$. In algebraic terms, this amounts to ask that no terms in the divergence-like part of $\det J_\Phi $ cancels with any term in its determinant-like part (maps with a linear part different from the identity can be reduced to the above form by multiplying by an invertible linear map).
Such a condition is satisfied under some symmetry conditions on $p(x,y)$ and $q(x,y)$, as in the case of $p(x,y)$ and $q(x,y)$ even polynomials, since the terms in $\psx + \qsy $ have odd degree, while the terms in $\psx \qsy - \psy \qsx  $ have even degree. 
In \cite{N} it was proved that such maps are globally invertible. A different proof for real maps, with some extensions, has been given in \cite{S}, by applying results proved in relation to the Global Asymptotic Stability Jacobian Conjecture \cite{F}, \cite{Gl}, \cite{Gu}.

In this paper we give a complete characterization of real divergence-free maps, proving that every such map has the form
\begin{equation}\label{divfree}
\Phi(x,y) = \left( x  + \sum_{i=2}^n \varepsilon_i \alpha (\beta x - \alpha y)^i, y + \sum_{i=2}^n \varepsilon_i \beta (\beta x - \alpha y)^i \right) , 
\end{equation}
with $\alpha, \beta  \in \R$, $\varepsilon_i \in \R$ for $i=2, \dots, n$. This means that the nonlinear part of such maps essentially depends on the single variable $\beta x - \alpha y$, so that after a change of variables the map has the form of a shear map, $(u,v) \mapsto (u,v-f(u))$. Our approach consists in solving the system of PDEs
\begin{equation}\label{PDE}
\psx + \qsy =0, \qquad   \psx \qsy - \psy \qsx = 0 ,
\end{equation}
in the class of homogeneous polynomials, then proving that only polynomials with nonlinear part as in  (\ref{divfree}) solve the system (\ref{PDE}). As a consequence, we prove that under some simple algebraic conditions, a jacobian map is the composition of a linear automorphism and a shear map.

\section{Results}

Let $\Phi : \R^2 \rightarrow \R^2$, $\Phi(x,y) = (P(x,y),Q(x,y))$ be a real polynomial map. Let $ J_\Phi$ be its jacobian matrix. We say that $\Phi $ is a {\it jacobian map} if its jacobian determinant $\det J_\Phi$ is a non-zero constant. Possibly replacing $\Phi(x,y)$ with $\Phi(x,y) - \Phi(0,0)$, which has the same jacobian determinant, we may assume that $\Phi(0,0) = (0,0)$. Moreover
one can compose $\Phi$ with the inverse matrix of $J_\Phi(0,0)$, obtaining another jacobian map with linear part coinciding with the identity. Let us call again $\Phi$ the map obtained by such operations.
Then we can write 
$$
P(x,y) = x + \sum_{i=2}^n p_i(x,y), \qquad Q(x,y) = y + \sum_{i=2}^n q_i(x,y), 
$$
where $p_i(x,y)$, $q_i(x,y)$ are homogeneous polynomials of degree $i$. 
Let us set $\overline \Phi (x,y) = \Phi(x,y) - (x,y)$. In other words, $\overline \Phi $ is the map consisting only of the nonlinear terms of $\Phi$. We say that  $\Phi(x,y)$ is a {\it divergence-free jacobian  map} if it satisfies the condition
$$
{\partial  \over \partial x}  \sum_{i=2}^n p_i(x,y) + {\partial  \over \partial y}   \sum_{i=2}^n q_i(x,y),  = 0.
$$
As observed in the introduction, if $\Phi$ satisfies the above condition, then $\det J_{\overline \Phi} = 0$. \\
\indent We say that  $\Phi$ is a {\it shear map} if its nonlinear terms are linear combinations of powers of a homogeneous polynomial of degree 1, as for the map $\Phi(x,y) = (3x-4y+(x-y)^2,-2x+y+(x-y)^2)$. The adjective {\it shear} comes from the fact the nonlinear part of that such maps substantially acts as a displacement along a given direction. The displacement's amount depends nonlinearly on the position of $(x,y)$.   \\
\indent  In next lemma we show that a homogeneous polynomial satisfying simultaneously the divergence and the determinant condition is essentially a function of a single variable. We start considering homogeneous polynomials.

\begin{lemma} \label{gradoj}
Let $p(x,y)$, $q(x,y)$ be homogeneous polynomials of the same degree $j > 1$, not both identically zero. If 
\begin{equation}\label{ipotesi}
\psx + \qsy = 0, \qquad \psx \qsy - \psy \qsx = 0,
\end{equation}
then there exist $\alpha, \beta, \zeta \in \R$, $\alpha ^2 + \beta ^2 \neq 0$, $\zeta \neq 0$, such that $p(x,y)= \zeta \alpha (\beta x - \alpha y)^j, q(x,y) = \zeta \beta (\beta x - \alpha y)^j$.
\end{lemma}
{\it Proof.}
Let us assume $q(x,y)$ not to be constant. Then there exists a line $l = \{(x,y) :  x = t \cos \theta , y = t \sin \theta\}$ such that $q(t \cos \theta , t \sin \theta )$ is non-constant. The gradient $\nabla q$ does not vanish at any point of $l$. By the second equality in (\ref{ipotesi}), in a neighbourhood of the point $(\cos \theta , \sin \theta )$ there exists an analytic function $\psi$ such that $p(x,y) = \psi(q(x,y))$.

 Then there exist $\alpha, \beta \in \R$, $\alpha^2 +  \beta^2 \neq 0$ such that
$$
\alpha t^j = p(t \cos \theta , t \sin \theta ) = \psi(q(t \cos \theta ,t \sin \theta ) ) =  \psi(\beta t^j).
$$
One has $\beta \neq 0$, otherwise $q(t \cos \theta ,t \sin \theta ) = \beta t^j$ would be constant. Setting $s = \beta t^j$, one has $t^j = \frac s\beta$, hence 
$$
\psi(s) =  \psi(\beta t^j) = \alpha t^j =  \frac{\alpha s}{\beta} .
$$
This proves that locally $\beta p(x,y) - \alpha q(x,y) = 0$. By the identity principle for polynomials, the same equality holds on all of $\R^2$. As a consequence, there exist a homogeneous polynomial $r(x,y)$ of degree $j$ and $k\in \R$, $k\neq 0$, such that $p(x,y) =k \alpha   r(x,y)$, $q(x,y) =k \beta  r(x,y)$.
Then one has 
$$
0 = \psx + \qsy = k \alpha \rsx + k \beta  \rsy = k (\alpha,\beta) \cdot \nabla r(x,y).
$$
This implies that the level curves of $r(x,y)$ are just the lines $\beta x - \alpha y = const.$, hence there exists a function $\rho$ such that $r(x,y) = \rho(\beta x - \alpha y)$. 
Choosing $(x,y) = \left(  \frac{\beta t}{\alpha^2 +\beta^2} , - \frac{\alpha t}{\alpha^2 +\beta^2} \right) $ one has $ \beta x - \alpha y = t$ and
$$
 \rho(  t) = r\left(  \frac{\beta t}{\alpha^2 +\beta^2} , - \frac{\alpha t}{\alpha^2 +\beta^2} \right) = h t^j.
$$
for some  $h \in \R$, $h \neq 0$.
In conclusion, $r(x,y) = \rho(\beta x - \alpha y) = h (\beta x - \alpha y)^j$. Setting $\zeta = hk$ gives the statement.
\hfill$\clubsuit$\bigskip

\begin{remark} \label{remark1} 
The above lemma shows that a divergence-free  jacobian map can be written as follows,
\begin{equation} \label{formaremark}
P(x,y) = x + \sum_{i=2}^n \zeta_i \alpha_i (\beta_i x - \alpha_i y)^i, \quad 
Q(x,y) = y + \sum_{i=2}^n \zeta_i \beta_i (\beta_i x - \alpha_i y)^i,  
\end{equation}
with  $\zeta_i \neq 0$, $\alpha_i^2 + \beta_i^2 \neq 0$ for $i=0, ... , n$. In fact, since \lq\lq divergence" terms do not cancel with \lq\lq  determinant" terms, every homogeneous couple $(p_j(x,y),q_j(x,y))$ satisfies the hypotehses of lemma  \ref{gradoj}.
\end{remark}

Something more can be said comparing terms with different degrees. Next lemma is concerned with couples of arbitrary functions both depending only on a first-degree homogeneous polynomial.

\begin{lemma} \label{duerette} 
Let us set $p(x,y) = a \phi(\beta x - \alpha y) + b \psi(\delta x - \gamma y)$, $q(x,y) = c \phi(\beta x - \alpha y) + d \psi(\delta x - \gamma y)$, with $\phi$, $\psi$  non-constant functions of class $C^1$. If (\ref{ipotesi}) holds, then the lines $\beta x - \alpha y=0$ and $\delta x - \gamma y=0$ coincide.
\end{lemma}
{\it Proof.}
One has 
$$
0 = \psx \qsy - \psy \qsx = (ad -bc) (\alpha \delta - \beta \gamma ) \phi' (\beta x - \alpha y) \psi' (\delta x - \gamma y).
$$
At every point where both $\phi' $ and $\psi' $ do not vanish one has  $(ad -bc) (\alpha \delta - \beta \gamma ) =0$. If $\alpha \delta - \beta \gamma = 0$, then the thesis is proved. If $ad-bc = 0$, then there exist $m, \eta, \omega \in \R$, $m\neq 0$, $\eta^2 + \omega^2 \neq 0$ and a function $r(x,y)$ such that $p(x,y) =m \eta   r(x,y)$, $q(x,y) =m \omega  r(x,y)$. Then, as in lemma  \ref{gradoj}, one has
$$
0 = \psx + \qsy = m \eta \rsx + m \omega  \rsy = m (\eta, \omega) \cdot \nabla r(x,y).
$$
As in lemma  \ref{gradoj}, this proves that  there exists one-variable function $\rho(t)$ such that $r(x,y) = \rho(\omega x - \eta y)$. Now, assume by absurd that the lines $\beta x - \alpha y=0$ and $\delta x - \gamma y=0$ do not coincide. Then one can make an invertible change of variables, $u=\beta x - \alpha y$, $v=\delta x - \gamma y$, and consider the functions of $u$ and $v$ obtained from $p(x,y), q(x,y), r(x,y)$ by applying the change of variables,
$$
p^*(u,v) = a \phi(u) + b \psi(v) = m \eta  \rho(k_1 u + k_2 v), 
$$
$$
q^*(u,v) = c \phi(u) + d \psi(v) = m \omega  \rho(k_1 u + k_2 v), 
$$
where $k_1, k_2 \in \R$, $k_1^2 +  k_2^2 \neq 0$. Differentiating both equalities with respect to $u$ one has
$$
a \phi'(u) = m \eta  k_1 \rho'(k_1 u + k_2 v), \qquad c \phi'(u)  = m \omega k_1  \rho'(k_1 u + k_2 v).
$$
If $k_2 \neq 0$, then the above equalities hold only if  $\phi'$ and $\rho'$ are constant, against the hypothesis, hence $k_2=0$. Similarly, differentiating  with respect to $v$ one proves that $k_1 = 0$. This leads to $k_1 = 0 = k_2$, contradiction.
\hfill$\clubsuit$\bigskip

\begin{theorem} \label{teorema} 
Let $\Phi(x,y) = (P(x,y),Q(x,y))$ be a divergence-free jacobian map of degree $n$. Then there exist $\alpha, \beta \in \R$, $\alpha^2 + \beta^2 \neq 0$, $\varepsilon_i \in \R$, such that
\begin{equation}\label{forma}
P(x,y) = x + \sum_{i=2}^n \varepsilon_i \alpha (\beta x - \alpha y)^i, \quad 
Q(x,y) = y + \sum_{i=2}^n \varepsilon_i \beta (\beta x - \alpha y)^i ,
\end{equation}
\end{theorem}
{\it Proof.}
Let $n$ be the degree of $\Phi$. Due to lemma  \ref{gradoj}, $P(x,y)$ and $Q(x,y)$ are as in remark \ref{remark1}.  Let us set  $ \alpha = \alpha_n$,  $\beta = \beta_n$, $\varepsilon_n = \zeta_n$. Since we assume $\Phi$ to be of degree $n$, $\alpha^2 + \beta^2 \neq 0$, $\varepsilon_n \neq 0$. \\
\indent 
We claim that for $i= 2, ... , n-1$, $\alpha_i = \alpha$, $\beta_i = \beta$. \\
\indent 
Let us proceed by decreasing degree. The functions $p(x,y) = \zeta_n \alpha (\beta x - \alpha y)^n +\zeta_{n-1} \alpha_{n-1} (\beta_{n-1} x - \alpha_{n-1} y)^{n-1}$, $q(x,y) =  \zeta_n \beta (\beta x - \alpha y)^n +  \zeta_{n-1} \beta_{n-1} (\beta_{n-1} x - \alpha_{n-1} y)^{n-1}$ satisfy the hypotheses of lemma \ref{duerette}, since their terms do not interact with those ones of lower degrees. 
Then, by lemma \ref{duerette}, the terms $\beta x - \alpha y$ and $\beta_{n-1} x - \alpha_{n-1} y$    are proportional. Since $\beta x - \alpha y$ does not vanish identically, there exist $\xi_{n-1} \in \R$ such that $\alpha_{n-1} = \xi_{n-1} \alpha $, $\beta_{n-1} = \xi_{n-1} \beta$. 
As a consequence, the terms of degree $\geq {n-1}$ in $P(x,y)$ and  $Q(x,y)$ are all function of $\beta x - \alpha y$. Setting $\varepsilon_{n-1} = \zeta_{n-1} \xi_{n-1}^n$, the form of the $(n-1)$-degree term in (\ref{formaremark}) is as in (\ref{forma}). \\
\indent 
Now we may repeat the procedure replacing $n-1$ with $n-2$. We apply again lemma \ref{duerette} to the functions 
$p(x,y) = \alpha \bigg( \varepsilon_n (\beta x - \alpha y)^n + \varepsilon_{n-1} (\beta x - \alpha y)^{n-1}\bigg) + \zeta_{n-2} \alpha_{n-2} (\beta_{n-2} x - \alpha_{n-2} y)^{n-2} $, 
$q(x,y) = \beta  \bigg( \varepsilon_n (\beta x - \alpha y)^n +  \varepsilon_{n-1}  (\beta  x - \alpha  y)^{n-1} \bigg)+  \zeta_{n-2} \beta_{n-2} (\beta_{n-2} x - \alpha_{n-2} y)^{n-2}$, in order to prove that $\beta x - \alpha y$ and $\beta_{n-2} x - \alpha_{n-2} y$ are proportional. Since $\beta x - \alpha y$ does not vanish identically, there  exist $\xi_{n-2} \in \R$ such that $\alpha_{n-2} = \xi_{n-2} \alpha $, $\beta_{n-2} = \xi_{n-2} \beta$. \\
\indent 
Such a procedure can be applied $n-2$ times in order to prove that every nonlinear term in $\Phi$ is actually a power of $\beta x - \alpha y$ multiplied by a constant $\varepsilon_i = \zeta_i \xi_i^{i+1}$, $i\neq n$.  \hfill$\clubsuit$\bigskip

In the above theorem it is not assumed that for $i=2, ..., n-1$ one has $\varepsilon_i \neq 0$. Actually some degrees may not appear in a divergence-free jacobian map, as in $\Phi(x,y) = (x - y^2 - y^5,y)$.

In \cite{N} and \cite{S} some simple algebraic conditions ensuring a jacobian map to be divergence-free were taken into account. We apply them to the present situation. 

Given a polynomial $P$, we write $d(P)$ for its degree, $o(P)$ for its order.  We say that a polynomial is {\it even} if it is the sum of even degree monomials,  {\it odd} if it is the sum of odd degree monomials. 
Similarly, we say that a polynomial is {\it x-even} if it contains only terms with even powers of $x$, {\it x-odd} if it contains only terms with odd powers of $x$.
We say that a non-negative integer is a {\it gap} of $P$ if it is the difference of the degrees of two distinct monomials in $P$. We denote by $G(P)$ the gap-set of $P$. As an example, the polynomial $P(x,y) = x^3 + y^3 + x^2y^2 + y^7$ has gap-set $G(P) = \{0, 1, 3, 4 \}$. If $P$ has exactly one monomial, or if $P$ is identically zero, then we say that it has empty gap-set. 

We say that the couple of polynomials $(P,Q)$ satisfies the {\it gap condition} if for every monomial $M$ in $P$, one has $d(M) -1 \not\in G(Q)$. The gap condition is not symmetric, as shown by the couple $(P,Q) = (x+y^2,x^6+y^2)$. In such a case one has $G(P) = \{ 1 \} $, $G(Q) = \{ 4 \} $, so that $(P,Q)$ satisfies the gap condition, but $(Q,P)$ does not.

We say that  $(P,Q)$ satisfies the {\it symmetric gap condition} if both  $(P,Q)$ and $(Q,P)$ satisfy the gap condition. 

\begin{corollary} \label{corollario1} 
Let $\Phi: \R^2 \rightarrow \R^2$ be a jacobian map of the type $\Phi(x,y) = (ax+by+p(x,y),cx+dy+q(x,y))$, $a,b,c,d \in \R$, $o(p) > 1$, $o(q) > 1$. If one of the following holds, \\
\indent i)   $\max \{d(p),d(q)\} <  o(p) + o(q) -1$,  \\
\indent  ii)  both $p(x,y)$ and $q(x,y)$ are even polynomials,  \\  
\indent iii) $p$ is odd, $q$ is even and $(p,q)$ satisfies the  gap condition,  \\
\indent iv) $(p,q)$ satisfies the symmetric gap condition,  \\
\noindent  then  $\Psi = \Phi \circ (J_\Phi(0,0))^{-1}$ is a shear map.
\end{corollary}
{\it Proof.}
Let us consider the map  $\Psi = \Phi \circ (J_\Phi(0,0))^{-1}$. One has $\Psi(x,y) = (x+p^*(x,y), y+q^*(x,y))$, with $o(p^*) > 1$, $o(q^*) > 1$. The composition with a linear map does not change the properties in $i), \dots, iv)$, hence $\Psi$ satisfies all of them.
In \cite{S} it was proved that any of the above conditions implies $\Psi$ to be a divergence-free map. Then the conclusion comes from theorem \ref{teorema}.
\hfill$\clubsuit$\bigskip

Next corollary is concerned with a kind of symmetry which is not preserved by composition with linear transformations.

\begin{corollary} \label{corollario2} 
Let $\Phi: \R^2 \rightarrow \R^2$ be a jacobian map of the type $\Phi(x,y) = (x+p(x,y), y+q(x,y))$, $o(p) > 1$, $o(q) > 1$. If one of the following holds, \\  
\indent i)  $p$ is $x$-even, $q$ is $x$-odd, \\  
\indent  ii) $p$ is $y$-odd, $q$ is $y$-even, \\
\noindent  then  $\Phi$ is a shear map.
\end{corollary} 
{\it Proof.}
By theorem 2 of \cite{S}, $\Phi$ is  a divergence-free jacobian map. Then one can apply theorem \ref{teorema}.
\hfill$\clubsuit$

\end{document}